\input{BoxedEPS.tex}
\SetTexturesEPSFSpecial
\HideDisplacementBoxes
\documentclass[12pt]{amsart}
\usepackage{amsmath,amsfonts,amssymb}

\newcommand{\C}{\mathbb C}

\newcommand{\Q}{\mathbb Q}
\newcommand{\Z}{\mathbb Z}
\newtheorem{theorem}{Theorem}[section]
\newtheorem{proposition}[theorem]{Proposition}
\newtheorem{lemma}[theorem]{Lemma}
\newtheorem{corollary}[theorem]{Corollary}

\newtheorem{question}{Question}

\def\proof{{\bf {\noindent}Proof. }}
\def\endproof{\hfill${\blacksquare}$\bigskip}
\def\K{\mathcal H} 

\def\bar{\overline}
\def\G{\Gamma}

\subjclass{20E36,20F32}
\keywords{exponential growth, trees, valuations, uniform exponential growth} 
\thanks{The authors thank the Mathematics Department at  Bielefeld University, Germany, for their
hospitality and support of this research project. G. Noskov supported by GIF-grant G-454-213.06/95, Gr-627/9}

\begin{document}
\thispagestyle{empty}

\renewcommand{\baselinestretch}{1.2}

\title[Uniform Growth]{\Large Uniform Growth,  Actions on Trees and $GL_2$}
\author[Alperin \& Noskov]{\large Roger C. Alperin and Guennadi A. Noskov}

\maketitle

\section{Exponential Growth}

 Choose a finite generating set $S=\{s_1,\cdots, s_p\}$ for the group $\G$;
define the $S-length$ of an element as $\lambda_S(g)=min\{\ n\ |\ g=s_1\cdots s_n, s_i\in
S\cup S^{-1}\}$. The growth function $\beta_n(S, \G)=|\{ g\ |\  \lambda_S(g)\le
n\}|$ depends on the chosen generating set. A group has
exponential growth if the growth rate, $\beta(S,\G)=limit_{n\rightarrow
\infty} \beta_n(S,\G)^{\frac{1}{n}}$ is strictly greater than 1.
In fact, for another finite generating set $T=\{t_1,\cdots, t_q\}$ for $\G$, if both
$max_j\lambda_S(t_j)\le L$ and $max_i\lambda_T(s_i) \le L$, then
$\beta_n(S,\G)\le \beta_{Ln}(T,\G)$ and also the symmetric inequality. It then
follows that $\beta(S,\G)^L\le\beta(T,\G)$ and $\beta(T,\G)^L\le\beta(S,\G)$.
Using these remarks, Milnor showed that  exponential growth is independent of the
generating set.

For a group with exponential growth we consider  $$\beta(\G)=inf_{S}
\beta(S,\G).$$ If $\beta(\G)>1$ then $\G$ is said to have {\it uniform
exponential growth}.

Gromov has asked if there is a group of exponential growth which is not of 
uniform exponential growth. Indications are that such a group will be hard to find. 
Recently, \cite{An, On}, it has been shown that all solvable groups which have exponential
growth are of uniform exponential growth.

\section{Generalities}

The following uses the fact that given a generating set $S$ for $\G$ the set
of elements of the subgroup $\K$, a subgroup of index $d$, which are words in
$S$ of length at most $2d-1$ give a generating set for $\K$.

\begin{proposition}[\cite{SW}] If a group $\G$ has a subgroup $\K$ of finite
index $d$, then $\beta(\G)\ge \beta(\K)^{\frac{1}{2d-1}}$. \label{SW}
\end{proposition}

The next proposition is elementary, but also a very useful result.
 
\begin{proposition} If a group $\G$ has homomorphic image which is of
uniform exponential growth then $\G$ has uniform exponential growth.
\end{proposition}

A group is called {\it large} if it has a homomorphism onto a free non-abelian group.
It is important to realize that  a free group of rank $n$ has 
uniform exponential growth of rate 
$\beta=2n-1$. 

\begin{corollary}
If the group  $\G$ is virtually large then $\G$ has uniform exponential growth.
\end{corollary}

There has been some interesting recent work on linear groups which are large, \cite{L}, \cite{MV}, in a sense  building on some
ideas from \cite{Se}.

We say that a group $\G$ has the {\it UF-property} (uniformly contains a free
nonabelian semigroup) if there is a constant $n_{\G} \ge 1$ such that
for every generating set $S$ of $\G $ there exist two elements (depending on $S$) in $\G$ of word length
$\le n_\G$ and freely generating a free semigroup of rank 2.

\begin{proposition}
If a group $\G$ has the UF-property then $\G$
is of uniform exponential growth.
\end{proposition}
\proof Let $S$ be an arbitrary finite generating set for $\G$ and let $S_0=\{g,h\}$ be the pair
of words of $S-length$ less than or equal to $n_\G$ and freely generating a free 
semigroup $\G _0$ of rank 2. 
We
have 
$$\beta(S,\G)^{n_\G }=\beta (S\cup S_0,\G )\ge \beta(S_0,\G _0)\ge 2$$ and the proof is complete. 
\endproof

\section{Action on Trees}

We use a variant on the usual ping-pong lemma to obtain  free semigroups.

\begin{lemma}[\cite{BH}]
\label{pp}
Let $\G_0$ be a group of isometries of a tree $X$ with the
generating set $\{g_1,\ g_2\}$ where $g_1$ $g_2$ are hyperbolic isometries with
distinct axes $A_1$, $A_2$. Then one of the four pairs $\{g_1^{\pm 1},g_2^{\pm 1}\}$
freely generate a free semigroup of rank two.
\end{lemma}

\proof If $A_1$, $A_2$ are disjoint then let $[a_1,a_2]$ be
the shortest geodesic joining them. Let $A_1^{+},A_2^{+}$ be the rays,
starting at $a_1,a_2,$ such that $g_1,\ g_2$ translate towards the ends of
these rays. Let $R_1^{+},\ R_2^{+}$ be the rays $A_1^{+},\ A_2^{+}$ with the
first unit length segments removed. Let $p_1,\ p_2$ be the geodesic projection
maps of $X$ onto $R_1^{+},\ R_2^{+}$ respectively. Set $X_i=p_i^{-1}(R_i^{+}),\ i=1,2.$ 
Clearly
$X_1\cap X_2=\emptyset$ and $g_i^n(X_1\cup X_2)\subset X_i,\ n\ge 1,i=1,2.$ We assert now that
$g_1,\ g_2$ freely generate a free semigroup of rank 2. Indeed let $w_1,\ w_2$ be (positive) words in
$g_1,\ g_2$ and suppose $w_1=w_2$ in $\G_0$. If the words start with the
same letter, the final segments are also equal in $\G_0$; hence by induction we get
that $w_1=w_2$ as formal words. Thus, we may assume that the words $w_1,\ w_2$
start with different letters, say $g_1,\ g_2$ respectively. Then the image of $X_1\cup X_2$ under
$w_1,\ w_2$ belongs to $X_1,\ X_2$ respectively; hence a contradiction.
\endproof

\begin{theorem} Suppose that $\G$ is a finitely generated subgroup of $\mathrm{GL}_2(K)$ 
over  the  field $K$. 
\begin{enumerate}
\item[(char $\ne 0$)] If $K$ has
nonzero characteristic and $\G$ has  exponential growth then $\G$  satisfies the
UF-property and consequently has  uniform exponential growth.
\item[(char $=0$)] If $K$ has
 characteristic zero and $\G$ has exponential growth then either $\G$  has  uniform exponential growth or $\G$ is
conjugate to a subgroup of $\mathrm{GL}_2(\mathcal{O})$, for $\mathcal{O}$ a ring of integers in an
algebraic number field.
\end{enumerate}
\end{theorem}

\proof Since $\G$ is finitely generated we may assume that $K$ is finitely generated. For
any discrete rank one valuation $v$ of $K$ let $X_v$ be the corresponding
Bruhat-Tits tree.  We consider the action of $\G$ (or a subgroup of index 2, acting without inversions) on each of the
Bruhat-Tits trees $X_v$ for each of these valuations and split the proof of the theorem into cases. 
\begin{enumerate}
\item For any valuation $v$, $\G$ has a fixed point on $X_v$ .  In this case, the ring $A$
generated by the traces of elements of $\G$ lies in every valuation ring $A_v$  of $K$. In
non-zero characteristic this ring is  finite, \cite{Se} I.6.2, generated over the prime field by roots of unity,
so in fact there are only finitely many traces and it follows that  $\G$ is finite, \cite{Wz} 1.20, (hence of eventually
constant growth), if the group acts irreducibly, or it is solvable if the group acts reducibly and hence also of uniform
exponential growth, \cite{On}. In characteristic zero, the ring of traces is the ring of algebraic integers
$\mathcal{O}$ in the algebraic closure of
$\Q$ in $K$; as in
\cite{Bs}, one can find a suitable module so that the group is conjugated into  $\mathrm{GL}_2(\mathcal{\bar{O}})$, for a
somewhat larger ring of integers
$\mathcal{\bar{O}}$.

\item 
\label{line} There is a valuation such that $\G$ has an invariant line on $X_v$.
Then by \cite{Se} II.1.3,  $\G$ is contained in a Cartan
subgroup, which is an extension of diagonal subgroup by cyclic of order 2.
Hence $\G$ is virtually abelian and has  polynomial growth.

\item  For some valuation $v$ there is neither a  fixed point nor an invariant line
on the Bruhat-Tits tree $X_v$.  We
prove that the image of $\G$  (that is modulo scalar matrices) satisfies the UF-property with a constant
$n_{\G}=6$ and hence so does the original group. 

Let $S$
be a finite generating set of $\G$. The subcases are as follows.
\begin{enumerate}
\item [(a)] $S$ contains a hyperbolic isometry, say $g$ and any $s\in S$ leaves the axis 
$A_g$ invariant. Then we are in case \ref{line}.

\item [(b)] $S$ contains a hyperbolic isometry, say $g$ and there is $s\in S$ such that $sA_g\neq A_g.$ The isometry
$h=sgs^{-1}$ is  hyperbolic  with the axis $sA_g\neq A_g,$ hence by Lemma \ref{pp} one of the four pairs $\{g^{\pm
1},h^{\pm 1}\}$ freely generate a free semigroup of rank 2 and the length of $g^{\pm
1}$ and $h^{\pm 1}$ is at most 3.

\item [(c)] $S$ does not contain a hyperbolic isometry, that is all isometries are
elliptic, thus any $s\in S$ has a fixed vertex in $X_v.$ Then there exists $g\in S\cup S^2$ which is hyperbolic,
\cite{Se}, I.6.4. Again we may assume there is $s\in S\cup S^2$ such that $sA_g\neq A_g.$
 The argument used for the previous case shows now that one of the four
pairs $\{g^{\pm 1},sg^{\pm 1}s^{-1}\}$ freely generate a free semigroup of
rank 2 and the length of $g^{\pm 1}$, respectively, $sg^{\pm 1}s^{-1}$ is at most 6 and the
proof is complete. \endproof
\end{enumerate}
\end{enumerate}

\section{Incidentals}

In the characteristic zero case, we have shown that the  finitely generated group $\G$ of exponential growth is of uniform
exponential growth or it is conjugate to a subgroup of
$\mathrm{GL}_2(\mathcal{O})$, for
$\mathcal{O}$ a ring of integers in an algebraic number field. This however is not a dichotomy. Certainly, there are many cases
of discrete subgroups of $\mathrm{GL}_2(\C)$ which are of uniform exponential growth; these arise as the
fundamental group of a hyperbolic manifold \cite{GN}, \cite{SW}.

For further analysis we  consider the case when every element of $\G$ is elliptic (considered as a group of Mobius
transformations).  By dint of a theorem of Lyndon-Ullman
\cite{LU} the group is conjugate (via stereographic projection) to a group of rotations of the
sphere.  However, for  a linear group which has compact closure,  every eigenvalue of every element must be on the unit circle.
If we are also in the case where $\G$ is a subgroup of $\mathrm{GL}_2(\mathcal{O})$, then under every embedding
$\sigma:\mathcal{O}\rightarrow \C$ the above property of eigenvalues holds for $\sigma(\G)$,
since ellipticity is preserved. Hence every eigenvalue is a root of unity by Dirchlet's Theorem. 
By passing to a subgroup of finite index (using compactness) we can assume the determinant of the matrices in this group are all
unity. Thus the trace of the matrices are in fact a sum of a root of unity and its complex conjugate. 
But now there are only finitely many traces since the fraction field of $\mathcal O$ has finite
dimension over $\Q$ and  only finitely many roots of unity can belong to a number field of given degree.
It now follows \cite{Alp}, that $\G$ is finite since we can pass to a subgroup of finite index missing these finitely many
possible traces $\ne 2$. But this subgroup of finite index has elements with all eigenvalues which are 1; thus every element is
unipotent and so must trivial by compactness.  

We summarize these remarks in the following statement.
\begin{corollary} Suppose that $\G$ is a finitely generated subgroup of $\mathrm{GL}_2(\C)$ of exponential growth and not
of uniform exponential growth then $\G$ can not consist entirely of elliptic elements.
\end{corollary}

\section{Related Open Problems}
\begin{question} Prove  uniform exponential growth for linear groups of exponential growth. 
For example  show  uniform exponential growth for  $SL_2(\mathcal O)$ where $\mathcal O$ is the ring of integers
of a  real quadratic number field.\end{question}
\begin{question} Prove that uniform exponential growth is a quasi-isometry invariant.\end{question}
\begin{question} Classify or otherwise characterize
the groups which have the property that every generating set has two elements which generate a free
semigroup. (if torsion-free are they of finite cohomological dimension?)
\end{question}

\bigskip
\noindent
Department of Mathematics and Computer Science\\ San Jose State University,
San Jose, CA 95192, USA\\
E-mail: alperin@mathcs.sjsu.edu
\medskip

\noindent
Institute of  Mathematics, Russian Academy of Sciences\\
Pevtsova 13,
Omsk 644099, RUSSIA\\
E-mail: noskov@private.omsk.su


\begin{thebibliography}{ABCDEF}

\bibitem{An}  R. C. Alperin, Uniform growth of polycyclic groups, preprint 2000

\bibitem{Alp}  R. C. Alperin, Two dimensional representations of groups with property FA, Proc. Amer. Math. Soc. 
108, No. 1, Jan. 1990, 283-4

\bibitem{Bs}  H. Bass, Finitely generated subgroups of GL$_2$, {\bf The
Smith Conjecture} edited by H. Bass and J. W. Morgan, Academic Press, 1984


\bibitem{BH} M. Bucher  and P. de la Harpe, Free products with amalgamation and HNN
extensions which are of uniformly exponential growth, preprint 2000
\bibitem{LU} R. C. Lyndon and J. L. Ullman, Groups of elliptic linear fractional transformations,   
 Proc. Amer. Math. Soc. 18, 1967, 1119-1124
\bibitem{GN} F. Grunewald and G.Noskov, Largeness of certain hyperbolic
lattices, preprint 2000
\bibitem{L} A. Lubotzky, Free quotients and the first betti number of some hyperbolic manifolds, Transformation Groups,
1,1996,71-82
\bibitem{MV} G. A. Margulis and E. B. Vinberg, Some linear groups virtually have a free quotient, J. Lie Theory,
10, 2000, 171-180 
\bibitem{On} D. Osin, The entropy of solvable groups, preprint 2000
\bibitem{SW} P. B. Shalen, P. Wagreich, Growth rates, $\Z_p$-homology, and
volumes of hyperbolic 3-manifolds,  Trans. Am. Math. Soc., 331, No. 2, 1992, 
895-917

\bibitem{Se} J-P. Serre, {\bf Trees}, Springer-Verlag, 1980.

\bibitem{Wz}B. A. F. Wehrfritz, {\bf Infinite Linear Groups}, Springer-Verlag, 1973.
\end{thebibliography}
\end{document}